\documentclass{article}
\usepackage[utf8]{inputenc}
\usepackage{amsmath,amssymb,hyperref}
\title{On Alphatrion's Conjecture about Hamiltonian paths in a hypercube}
\author{Steppan Konoplev}
\date{}

\begin{document}
\maketitle

\section{Introduction}
Alphatrion conjectured [1] that it is possible to label the vertices of an $n$-dimensional hypercube with distinct positive integers such that for every Hamiltonian path $a_1, \dots, a_{2^n},$ we have $a_i + a_{i+1}$ prime for all $i.$ For a labeling of the vertices $V$ of a graph by a function $L: V \to \mathbb{N},$ let the edge sum of an edge $\overline{uv}$ be $L(u)+L(v).$ We prove that a graph $G = (V, E)$ can be labeled with distinct positive integers such that the edge sum for all $e \in E$ is prime if and only if $G$ is bipartite. Since the hypercube graph $Q_n$ is embeddable in the bipartite graph $K_{2^{n-1}, 2^{n-1}},$ this settles Alphatrion's conjecture in the affirmative.
\section{Proof of General Theorem}
Any graph not embeddable into $K_{m, n}$ for some $m, n$ is not bipartite, hence it contains an odd cycle [2]. Suppose the vertices in that cycle are labeled $a_1, \dots, a_k.$ Then $2(a_1+\dots+a_k) = (a_1+a_2)+\dots+(a_k+a_1)$ is both even and the sum of an odd number of primes. Thus, one of the primes must be $2,$ which is only possible via $2=1+1.$ But the vertex labels are distinct, so this cannot happen.\\\\
To label $K_{m,n},$ let $a, a+d, \dots, a+Ld$ be an arithmetic progression of primes with $L \ge 2mn-2m-n+3;$ such a progression is guaranteed to exist by the Green-Tao Theorem [3]. Represent $K_{m,n}$ as $A \cup B$ where $A = \{a_0, \dots, a_{m-1}\}, B = \{b_0, \dots, b_{n-1}\},$ and there is an edge between $a_i$ and $b_j$ for all $i,j.$ Let $c = a+(mn-m-n+2)d$ and let $a_k = (k(n-1)+1)d, b_k = c+kd.$ Then $a_i+b_j = a+(mn-m-n+3+j+i(n-1))d$ and $mn-m-n+3+j+i(n-1) \le mn-m-n+3+j+i(n-1) \le mn-m-n+3+(n-1)+(m-1)(n-1) = 2mn-2m-n+3 \le L,$ so all edge sums are prime. Furthermore, $a_i \le a_{m-1} = (mn-m-n+2)d < c = b_0 \le b_j,$ so all edge labels are distinct.
\section{References}

1. \url{https://artofproblemsolving.com/community/c907967h1879578_hamilton_paths__primes_} Image of the post available at \url{https://imgur.com/a/8XiUfbj} for those without an account.\\
2. The proof is due to Kőnig. If $G = A \cup B$ is bipartite ($A, B$ are independent sets), any path starting at $a \in A$ bounces back and forth between $A$ and $B,$ so it can only return to $A$ and hence $a$ after an even number of steps, hence any cycle must have even length; similarly for any path starting at $b \in B.$ Suppose $G = A \cup B$ is bipartite and has an odd cycle $a_1 \to a_2 \to \dots \to a_k \to a_1.$ WLOG $a_1 \in A.$ Then $a_2 \in B,$ so $a_3 \in A,$ and eventually $a_k \in A.$ But $a_k$ is connected to $a_1,$ contradiction. Thus, a graph is bipartite iff it has no odd cycles.\\
3. Green, Ben, and Terence Tao. "The primes contain arbitrarily long arithmetic progressions." Ann. Math (2005). 
\end{document}